\documentclass[a4paper,11pt]{article}

\usepackage[a4paper,margin=2cm]{geometry}

\usepackage[english]{babel}

\usepackage{amsfonts,amssymb}
\usepackage{times}

\usepackage{amsthm}
\newtheorem{theorem}{Theorem}
\newtheorem{proposition}[theorem]{Proposition}
\theoremstyle{definition}
\newtheorem{remark}[theorem]{Remark}

\usepackage{mathtools}

\usepackage{graphicx}

\usepackage{doi}
\usepackage{hyperref}
\hypersetup{
    colorlinks=true,
    linkcolor=blue,      
    urlcolor=cyan,
    citecolor=magenta,
    }

\usepackage{xcolor}

\usepackage{authblk}

\usepackage[nocompress]{cite}

\newcommand{\keywords}[1]{%
  \begin{quote}
    \noindent\small\textbf{Keywords:} #1
  \end{quote}
}

\newcommand{\msc}[1]{%
  \begin{quote}
    \noindent\small\textbf{MSC2020 Classification:} #1
  \end{quote}
}

\begin{document}

\title{Yet another note on notation}
\author{Mircea Dan Rus}
\affil{Technical University of Cluj-Napoca, Department of Mathematics\\\texttt{rus.mircea@math.utcluj.ro}}
\date{}

\maketitle

\begin{center}\small\itshape This article has been accepted for publication in the
International Journal of Mathematical Education in Science and Technology,
published by Taylor \& Francis.\end{center}
%% After publication, replace the statement above with the Version-of-Record notice:
%% \begin{center}\small\itshape This is an original manuscript of an article published by
%% Taylor \& Francis in the International Journal of Mathematical Education in Science and
%% Technology on [date of publication], available online:
%% \url{https://doi.org/10.1080/0020739X.2026.2701739}.\end{center}

\begin{abstract}
Back in 1755, Euler explored an interesting array of numbers that now frequently appears in polynomial identities,
combinatorial problems, and finite calculus, among other places. These numbers share a strong connection with well-known
number families, such as those of Stirling, Bernoulli, and Fubini. Despite their importance, they often go unnoticed
because of the lack of a specific name and standard notation. This paper aims to address this oversight by proposing
an appropriate name and notation, aligned with established mathematical conventions, and supported by (we hope)
strong enough arguments to facilitate their acceptance from the mathematical community.
\end{abstract}


\msc{Primary 05-01; Secondary 05A10, 00A35.}

%% MSC 2020:
%% 05-01   	Introductory exposition (textbooks, tutorial papers, etc.) pertaining to combinatorics
%% 05A10   	Factorials, binomial coefficients, combinatorial functions
%% 00A35   	Methodology of mathematics

\section{Introduction}

An early milestone in the study of combinatorics is realising that not
everything can be counted using elementary principles or expressed with basic
formulas. For example, when counting surjective functions that map an
$m$-set (that is, a set with $m$ elements) to an $n$-set, we come across the
(in)famous formula $\sum_{k=1}^{n}(-1)^{n-k}\binom{n}{k}k^{m}$. To the
uninitiated, this may come as a surprise, given that counting injective
functions --- a seemingly similar problem --- results in a straightforward and
simple formula.

When I discuss this problem with my students, typically in relation to the
inclusion--exclusion principle, their reactions vary from joy and relief at
finding a solution,\footnote{See, for example, \cite[Ch.
3]{Mazur2010}.} to dissatisfaction with the method or the formula itself.
Restating the problem in more familiar terms, such as counting the ways of distributing
$m$ distinct toys to $n$ children with each child receiving at
least one toy,\footnote{The correspondence between each toy and the child who
receives it represents a surjective function.} helps a little, but not enough.
To ease the frustration of students, I decided a few years ago to
\textquotedblleft pack\textquotedblright\ the formula into a notation and
\textquotedblleft unpack\textquotedblright\ it only when needed, such as
during computations or other algebraic manipulations that required the actual expression.

Those familiar with the calculus of finite differences may have already
recognised the above formula as the so-called \emph{differences of zero},
represented as $\Delta^{n}x^{m}(0)$, or $\Delta^{n}0^{m}$ for short. However, this was not
what I had in mind, as the context was different. I needed a notation that
captured the combinatorial nature of the formula in connection with similar
problems whose results were already represented with standard familiar
notations. My initial research yielded no existing notation for this purpose,
so I decided to create my own notation, and then presented it to my students,
along with some supporting arguments.

Upon further examination, I discovered that the numbers expressed by the formula above
have a long history and have been independently rediscovered in
various contexts; the first seems to have been Euler \cite[\S 14]{Euler1755} in 1755.

As will be shown next, the use of \emph{the notation}
simplifies many results involving these numbers while also shedding new light
on their significance, potentially revealing deeper insights.

\section{The context}

Combinatorial counting often results in expressions that include products,
particularly powers. So, let us start from there.

According to Cajori \cite[\S 298--309]{Cajori1993}, the modern form of
exponential notation, although initially restricted to positive integral
exponents, dates back to Descartes in his \emph{G\'{e}om\'{e}trie}~\cite{Descartes1637} (1637).
This was followed by the introduction of negative and fractional exponents by
Wallis~\cite{Wallis1656} (1656) and Newton~\cite{Newton1676} (1676). In 1743, Euler's most celebrated formula
$\mathrm{e}^{ix}=\cos x+i\sin x$ showcased the power of this notation, marking
the first appearance in print of imaginary\footnote{The use of the symbol $i$
for the imaginary unit is also attributed to Euler, but it appeared in print
only in 1794, after his death (see \cite[\S 498]{Cajori1993}).} exponents. As
Cajori noted, \textquotedblleft there is perhaps no symbolism in ordinary
algebra which has been as well-chosen and is as elastic as the Cartesian
exponents\textquotedblright.

Besides the usual integral powers, there exist also \emph{factorial powers}: the
\emph{rising factorial }$x^{\overline{n}}=x(x+1)\cdots(x+n-1)$ that gets this
notation from Capelli \cite{Capelli1893} in 1893, and the \emph{falling
factorial }$x^{\underline{n}}=x(x-1)\cdots(x-n+1)$ written here in the
notation of Toscano \cite{Toscano1939} from 1939. These notations are
mnemonic: the exponent $n$ represents the number of factors, while the bar
above or below $n$ indicates that the factors increase or decrease stepwise
from $x$, respectively. Graham, Knuth, and Patashnik \cite[p. 48]{Graham1994}
anticipated three decades ago that \textquotedblleft the underline/overline
convention is catching on, because it's easy to write, easy to remember, and
free of redundant parentheses\textquotedblright.\footnote{Non-mnemonic
notations in use for $x^{\underline{n}}$ include $x_{(n)}$, $x^{(n)}$,
$x_{\left\langle n\right\rangle }$ or $[x]_{n}$. Furthermore,
\textquotedblleft Pochhammer's symbol\textquotedblright\ $(x)_{n}$ is used to
represent both $x^{\underline{n}}$ and $x^{\overline{n}}$ in different
contexts (see \cite{Knuth1992} for more details).}

Factorial powers are ubiquitous in finite calculus, combinatorics, and the
theory of special functions. They also share characteristics with the more
common integral powers, as presented in Vandermonde's\footnote{Vandermonde used
$[x]^{n}$ to denote the falling factorial power $x^{\underline{n}}$.}
\emph{M\'{e}moire} \cite{Vandermonde1772} from 1772. For instance, the
binomial theorem still holds if we consistently replace integral powers
with one of factorial powers (see \cite[p. 71, Exercise 33]{Knuth1997}).
In addition, in the calculus of finite differences, the falling factorial
$x^{\underline{n}}$ plays the same role as the monomial $x^{n}$ in the
differential calculus.

Revisiting the combinatorial context, the power $m^{n}$ counts the
permutations of any $n$ selections from a given set of $m$ distinctively
labelled objects, allowing the same object to be selected more than once. In
contrast, the falling factorial $m^{\underline{n}}$ provides the answer when
the selections must be distinct, which means that no repetitions are allowed. This
corresponds to counting functions and injective functions, respectively, that
map an $n$-set to an $m$-set, since each function can be identified with the
sequence of its values (provided that a total order on the domain of the functions is specified).
In this context, it is natural to ask: why not have
\textbf{an exponential notation also for the number of surjective functions}
as well?

We need to clarify one point before proceeding. If $m<n$, then the number of
injective functions from an $n$-set to an $m$-set is clearly $0$, which is
also the value of $m^{\underline{n}}$, since $0=m-m$ is one of the factors.
The formula $m^{\underline{n}}$ generally holds, but is only of real interest
when $m\geq n$. Under this condition, it makes sense to count the surjective
functions from an $m$-set to an $n$-set, with the roles of the sets reversed
from the injective case. Thus, to maintain consistency with the factorial
power notation and the analogy with the number of injective functions, $m$ and
$n$ must switch places.

\section{The notation}

The proposed notation for the number of surjective functions from an $m$-set
to an $n$-set is $n^{\{m\}}$. The power notation is motivated by its correspondence
with the number $n^{m}$ of all functions between the same sets and by the
actual expression
\begin{equation}
n^{\{m\}}=\sum_{k=0}^{n}(-1)^{n-k}\binom{n}{k}k^{m} \label{yanon_eq:01}%
\end{equation}
which is a linear combination\footnote{The coefficients are independent of the
exponent $m$.} of the powers $0^{m}$, $1^{m}$, $2^{m}$, \ldots, $n^{m}$. The
reasons for including also the powers of $0$ will soon become clear, as we
consider (\ref{yanon_eq:01}) even when $m$, or $n$, or both are $0$.

The braces around the exponent could suggest that we use the
inclusion--exclusion principle (which has to do with intersections and reunion
of sets) to obtain the formula. But, mainly, the braces are replicated from
the notation $\genfrac{\{}{\}}{0pt}{}{m}{n}$ for a related set of numbers, called \emph{the
Stirling set numbers} (also commonly known as Stirling numbers of the second
kind)\footnote{There exist also \emph{the Stirling cycle numbers / of the first kind},
that count the number of ways to partition a finite set of a given size into a specified number of cycles.}
which count the number of ways to partition an $m$-set into $n$ parts.
These numbers were first considered in a purely algebraic context by Stirling
in 1730 in his influential book \emph{Methodus Differentialis}
\cite{Stirling1730}, in the expression of monomials as linear combinations of
falling factorials:%
\begin{equation}
x^{m}=\sum_{n=0}^{m}\genfrac{\{}{\}}{0pt}{}{m}{n}x^{\underline{n}}, \label{yanon_eq:02}%
\end{equation}
only to be rediscovered in the combinatorial context of set partitions by Saka
in 1782 (according to Knuth \cite[p. 504]{Knuth2011}).

In analogy to the expression of the binomial coefficients $\binom{m}{n}$ in
terms of the falling factorial powers $m^{\underline{n}}$%
\begin{equation}
\binom{m}{n}=\frac{m^{\underline{n}}}{n!}, \label{yanon_eq:03}%
\end{equation}
the Stirling set numbers can be expressed as
\begin{equation}
\genfrac{\{}{\}}{0pt}{}{m}{n}=\frac{n^{\{m\}}}{n!}, \label{yanon_eq:04}%
\end{equation}
justifying the choice for connecting the two notations with the use of braces.
Note that (\ref{yanon_eq:04}) is an immediate consequence of a two-step construction
of any surjective function. First, partition the $m$-set into $n$ parts,
which can be done in $\genfrac{\{}{\}}{0pt}{}{m}{n}$ ways. Then, assign a different output value
from the $n$-set to each part, which can be done in $n!$ ways. Consequently,
$n^{\{m\}}=n!\,\genfrac{\{}{\}}{0pt}{}{m}{n}$ (see also \cite[p. 75]{Stanley2012}).

\begin{remark}
The notation for the Stirling set numbers can be traced back to
Karamata \cite{Karamata1935} in 1935 and was further promoted
by Knuth \cite{Knuth1992}\footnote{Yet, the acknowledgment
of Karamata's primacy came later in \cite{Graham1994}.} who remarked that
\textquotedblleft the meaning of
$\genfrac{\{}{\}}{0pt}{}{m}{n}$ is easily remembered, because braces $\{~\}$ are commonly used to denote
sets and subsets\textquotedblright. Knuth also pointed out that the
\textquotedblleft lack of a widely accepted way to refer to these numbers has
become almost scandalous\textquotedblright. Unfortunately, after three decades,
there has been little progress toward reaching a consensus.
\end{remark}

You may wonder why introduce a new notation when we can use
$n!\,\genfrac{\{}{\}}{0pt}{}{m}{n}$, which is combinatorial in nature and fits the context very well?
There are several reasons why a proper notation is needed for these
\textquotedblleft multiples of the Stirling set numbers\textquotedblright.
First, there is an obvious precedent. Falling factorial powers and binomial
coefficients each have their own distinct representation. However, the connection
between them expressed in (\ref{yanon_eq:03}) is identical to that in
(\ref{yanon_eq:04}). Second, \textquotedblleft$n!\,\genfrac{\{}{\}}{0pt}{}{m}{n}$\textquotedblright\ is
just another formula, not a notation, which forces us to write and think
specifically in terms of Stirling's numbers. For a more contrasting picture,
imagine mathematics without a notation for the cosine (or even
without the function itself), simply because \textquotedblleft
notations\textquotedblright\ such as $\sin\left(  \frac{\pi}{2}-x\right)  $ or
$\sin\left(  \frac{\pi}{2}+x\right)  $ can be used instead. Is $\sin$ not
enough\footnote{No pun intended.}?

To further argue that the proposed notation is a benefit rather than an
overcomplication, we need to delve deeper.

\section{The name}

\textquotedblleft A new symbol in algebra is only half a benefit unless it
[has] a new name\textquotedblright, wrote Whitworth \cite{Whitworth1877} who
in 1877 proposed the name \textquotedblleft sub-factorial $n$%
\textquotedblright\footnote{Nowadays, it is also spelled \emph{subfactorial}%
.}\ for the expression $n!\sum_{k=0}^{n}\frac{(-1)^{k}}{k!}$. This represents
the number of \emph{derangements} of $n$ ordered objects, i.e., permutations
where no object remains in its original position. The name stuck, while
Whitworth's notation\footnote{This was a derivation from $\lfloor
\underline{n}$ which was one of the notations for the factorial of $n$.} which
looked like $\underline{\Vert n}$ was in time replaced by $!n$ and
$n\mathord{\text{\textexclamdown}}$, each meant as a \emph{symbolic
derangement} of \textquotedblleft$n!$\textquotedblright.

To propose a name for the numbers $n^{\{m\}}$, let us revisit the subfactorial
formula and rewrite it as
\begin{equation}
!n=\sum_{k=0}^{n}(-1)^{n-k}\binom{n}{k}k! \label{yanon_eq:05}%
\end{equation}
which has a structure identical to that of $n^{\{m\}}$ in (\ref{yanon_eq:01}),
except that the powers are replaced by factorials. In the spirit of
\textquotedblleft Whitworth's naming principle\textquotedblright\ of using the
prefix \textquotedblleft\emph{sub-\textquotedblright} to indicate that the
subfactorial is subordinated to the factorial, the numbers $n^{\{m\}}$ could
be named \emph{sub-powers }(with or without a hyphen).

We could also refer to these numbers as the \emph{ordered Stirling set
numbers}, which reflects the combinatorial interpretation of the formula
$n^{\{m\}}=n!\,\genfrac{\{}{\}}{0pt}{}{m}{n}$. To underline that we use a power-like notation,
we will go with \textquotedblleft subpowers\textquotedblright(which is also shorter!).

We see the \emph{subordination} relation explicitly reflected in the formulas%
\begin{align}
n!  &  =\sum_{k=0}^{n}\binom{n}{k}~!k\label{yanon_eq:06}\\
n^{m}  &  =\sum_{k=0}^{n}\binom{n}{k}k^{\{m\}} \label{yanon_eq:07}%
\end{align}
which express the factorial as a combination of subfactorials and the usual
power as a combination of subpowers. These relations can be derived from some
simple combinatorial arguments, independent of (\ref{yanon_eq:05}) and
(\ref{yanon_eq:01}). Both sides of (\ref{yanon_eq:06}) count the permutations of $n$
ordered objects, where the sum is obtained by grouping the permutations having
the same number ($k$) of objects not in their initial position. Similarly,
(\ref{yanon_eq:07}) counts the number of functions from an $m$-set to an $n$-set
in two ways, where the term $\binom{n}{k}k^{\{m\}}$ counts the functions that
take exactly $k$ distinct values.

\begin{remark}
\label{rem:sumation_range}The range of $k$ that corresponds to the nonzero
terms in (\ref{yanon_eq:07}) is from $1$ to $\min\{m,n\}$, if $m\neq0$. For $m=0$,
there is only one term $\binom{n}{0}0^{\{0\}}$. In some cases, it is more
convenient to consider the summation in a range independent of $n$
(e.g., from $0$ to $m$, or from $1$ to $m$), since $\binom{n}{k}=0$ for $k>n$.
\end{remark}

It is worth mentioning that once (\ref{yanon_eq:06}) and
(\ref{yanon_eq:07}) are established, the results in (\ref{yanon_eq:05}) and (\ref{yanon_eq:01})
follow as immediate consequences of a more general property, usually called
\emph{binomial inversion} (see \cite[p. 198]{Stanley2012}), to which we will refer
several times in this paper.

\begin{proposition}
\label{yanon_binom_inv}
Given two sequences\footnote{The sequences may also
be finite, but they must have the same number of terms.} $(a_{n})_{n\geq0}$ and
$(b_{n})_{n\geq0}$, the following statements are equivalent:

\noindent{\upshape(i)}$\quad\displaystyle b_{n}=\sum_{k=0}^{n}\dbinom{n}{k}a_{k}$,\quad
for all $n$;\hfill{\upshape(ii)}$\quad\displaystyle a_{n}=\sum_{k=0}^{n}(-1)^{n-k}%
\dbinom{n}{k}b_{k}$,\quad for all $n$.
\end{proposition}

In the terminology used in \emph{The On-Line Encyclopedia of Integer
Sequences} (OEIS) \cite{OEIS}, $(b_{n})_{n\geq0}$ is called the \emph{binomial
transform} of $(a_{n})_{n\geq0}$, while $(a_{n})_{n\geq0}$ is the
\emph{inverse binomial transform} of $(b_{n})_{n\geq0}$.

In consequence, for any fixed exponent, the resulting sequence of subpowers
can be seen as the inverse binomial transform of the corresponding sequence of powers.

\section{The first impression}

To get a first taste of the subpower concept and to break the initial
resistance to the introduction of a new notation, let us explore some basic
properties of these numbers in light of the new symbolism and terminology.

\subsection{The values}

When $n$ is small, values such as%
\[
0^{\{m\}}=0,\quad1^{\{m\}}=1,\quad2^{\{m\}}=2^{m}-2,\quad3^{\{m\}}%
=3^{m}-3\cdot2^{m}+3\qquad(m\geq1)
\]
are straightforward manifestations of the general formula, although they can
be derived independently. In contrast, when $n$ is large (compared to $m$), it
follows from the combinatorial definition that%
\begin{equation}
n^{\{m\}}=\left\{
\begin{tabular}
[c]{rl}%
$0\,$, & $\text{if }n>m$\\
$n!\,$, & $\text{if }n=m$%
\end{tabular}
\right.  \qquad(m,n\geq0). \label{yanon_eq:08}%
\end{equation}
In particular, $n^{\{0\}}=\delta_{0,n}$ if using the Kronecker delta function.
The value $0^{\{0\}}=1$ agrees also with the combinatorial and algebraic
interpretation\footnote{For a brief account about the \textquotedblleft$0^{0}$
controversy\textquotedblright, we point to \cite[p.
407--408]{Knuth1992}.} that $0^{0}=1$.

The result in (\ref{yanon_eq:08}) with $n^{\{m\}}$ written explicitly is well known
and can be attributed to Euler \cite[p. 31--33]{Euler1769} in a
noncombinatorial context. A detailed account given by Gould on the
implications of this result can be found in \cite{Gould1978}.

Extending (\ref{yanon_eq:08}) further to $m=n+1$ or $m=n+2$ produces the
expressions
\begin{align}
n^{\{n+1\}}  &  =n!\dbinom{n+1}{2}=\dfrac{(n+1)!\cdot n}{2}\label{yanon_eq:09}\\
n^{\{n+2\}}  &  =n!\left(  \dbinom{n+2}{3}+3\dbinom{n+2}{4}\right)
=\dfrac{(n+2)!\cdot n(3n+1)}{24} \label{yanon_eq:10}%
\end{align}
through a direct combinatorial approach.

To prove (\ref{yanon_eq:09}), consider the problem of counting in how many ways
$n+1$ toys can be divided among $n$ children, ensuring that each child is happy
(that is, everyone receives at least one toy). In each configuration, there is
exactly one child (the happiest!) who receives two toys, while the others each
get just one toy. Choosing the happiest child can be done in $n$ ways, while
his toys can be chosen in $\binom{n+1}{2}$ ways. The distribution of the remaining
$n-1$ toys to $n-1$ children can be done in $(n-1)!$ ways. The product of
these values gives the expression for $n^{\{n+1\}}$.

The result in (\ref{yanon_eq:10}) is obtained similarly, with $n+2$ toys distributed
among $n$ children such that all are happy, focusing on the child or
children who receive more than one toy. The details are left as an exercise to
the interested reader.

\subsection{The recurrence}

In analogy to the usual powers and the falling factorial powers that verify
the recurrences
\[
n^{m}=n\cdot n^{m-1},\qquad n^{\underline{m}}=n\cdot(n-1)^{\underline{m-1}}~\qquad(m,n\geq1),
\]
the subpowers behave like a \textquotedblleft superposition\textquotedblright of the previous two, satisfying%
\begin{equation}
n^{\{m\}}=n\left(  n^{\{m-1\}}+(n-1)^{\{m-1\}}\right)  \qquad(m,n\geq1).
\label{yanon_eq:11}%
\end{equation}
There is a simple combinatorial argument for (\ref{yanon_eq:11}) that uses the
previous model of distributing toys to make children happy. With $m$ toys and
$n$ children, let us distinguish a \emph{special} toy (e.g., the most
expensive one). There are $n$ ways to select the child who receives it, and we
will call this child \emph{lucky}. Distributing the remaining $m-1$ toys
depends on whether the lucky child is happy with just one toy or not. When he
is not happy, then we must make all the $n$ children happy using $m-1$ toys,
leading to $n^{\{m-1\}}$ ways of distributing the remaining toys among them.
If the lucky child is happy with just the special toy, then he receives no
other toy, while the remaining $m-1$ toys can be distributed in
$(n-1)^{\{m-1\}}$ ways among the other $n-1$ children, ensuring all are happy.
As a result, there are $n$ ways of giving the special toy, then $n^{\{m-1\}}%
+(n-1)^{\{m-1\}}$ ways of distributing the remaining toys, which leads
to (\ref{yanon_eq:11}).

The immediate practicality of (\ref{yanon_eq:11}) is evident in the recursive
generation of the array $\left(  n^{\{m\}}\right)  _{m,n\geq0}$, starting from
the initial values $0^{\{0\}}=1$ and $0^{\{m\}}=n^{\{0\}}=0$ for $m,n\geq1$.
This array of numbers read by rows (with $m$ as the row index) and omitting
the (zero) values for $n>m$ appears as the entry A131689 in the OEIS
\cite{OEIS}.

\begin{table}[h]
\caption{The first few lines of the array $\left(  n^{\{m\}}\right)
_{m,n\geq0}$. The zero entries are omitted.}%
\setlength{\tabcolsep}{6pt}
\renewcommand{\arraystretch}{1.2}
\centering
\par%
\begin{tabular}
[c]{c|rrrrrrrrr}%
\multicolumn{1}{l|}{\hspace*{-4pt}\raisebox{-2pt}{$m$}\hspace{-6pt}%
{\Large $\diagdown$}\hspace{-6pt}\raisebox{4pt}{$n$}\hspace*{-5pt}} &
\multicolumn{1}{c}{{\small \textbf{0}}} &
\multicolumn{1}{c}{{\small \textbf{1}}} &
\multicolumn{1}{c}{{\small \textbf{2}}} &
\multicolumn{1}{c}{{\small \textbf{3}}} &
\multicolumn{1}{c}{{\small \textbf{4}}} &
\multicolumn{1}{c}{{\small \textbf{5}}} &
\multicolumn{1}{c}{{\small \textbf{6}}} &
\multicolumn{1}{c}{{\small \textbf{7}}} &
\multicolumn{1}{c}{{\small \textbf{8}}}\\\hline
{\small \textbf{0}} & $1$ &  &  &  &  &  &  &  & \\
{\small \textbf{1}} &  & $1$ &  &  &  &  &  &  & \\
{\small \textbf{2}} &  & $1$ & $2$ &  &  &  &  &  & \\
{\small \textbf{3}} &  & $1$ & $6$ & $6$ &  &  &  &  & \\
{\small \textbf{4}} &  & $1$ & $14$ & $36$ & $24$ &  &  &  & \\
{\small \textbf{5}} &  & $1$ & $30$ & $150$ & $240$ & $120$ &  &  & \\
{\small \textbf{6}} &  & $1$ & $62$ & $540$ & $1560$ & $1800$ & $720$ &  & \\
{\small \textbf{7}} &  & $1$ & $126$ & $1806$ & $8400$ & $16800$ & $15120$ & $5040$ & \\
{\small \textbf{8}} &  & $1$ & $254$ & $5796$ & $40824$ & $126000$ & $191520$ & $141120$ & $40320$%
\end{tabular}
\end{table}

The first to have discovered and investigated this array of numbers and the recurrence
that governs it --- though unaware of its combinatorial significance
or its connection with the Stirling set numbers ---
appears to have been Euler. In 1755, Euler \cite[\S 14]{Euler1755} tabulated
these numbers, up to $n=7$ and $m=8$, and used them to express the forward
differences of the monomial $x^{m}$ in terms of the other monomials (as
will be seen in (\ref{yanon_eq:17})). This property of the subpowers will be
developed in the next section.

\subsection{The binomial expansion}

A power-like notation must exhibit power-like properties. The ultimate test is
the binomial expansion. Surprisingly, the formula holds:%
\begin{equation}
(a+b)^{\{m\}}=\sum_{k=0}^{m}\binom{m}{k}a^{\{k\}}b^{\{m-k\}}\qquad
(a,b,m\geq0), \label{yanon_eq:12}%
\end{equation}
which is a decisive argument in support of a power-like notation. Once more,
the combinatorial interpretation of the two sides of the formula
provides an immediate proof. Consider the problem of counting in how many ways
we can make $a$ girls and $b$ boys happy by distributing $m$ toys among them.
Clearly, the direct answer is the left-hand side of (\ref{yanon_eq:12}), while the
right-hand side provides the answer by grouping the configurations based on
the number $k$ of toys received by the girls.

\begin{remark}
The actual range of the index $k$ in (\ref{yanon_eq:12}), corresponding to the
nonzero terms, is from $a$ to $m-b$. However, the full range from $0$ to $m$
is preferred for its analogy with the binomial formula for positive integral powers.
\end{remark}

\section{The bigger picture. Making connections}

The concept of subpowers would be of limited use if it were confined solely to
expressing the results of certain combinatorial problems. In what follows, we
will see that the subpower notation fits quite naturally in several other contexts.

\subsection{Finite differences}

As already mentioned, Euler \cite[\S 13--15]{Euler1755} encountered the
subpowers in the computation of finite differences for polynomials. First, we
need to give a brief introduction to this subject.

For any nonzero number $h$, let $\Delta_{h}$ be \emph{the forward difference
operator} that maps a function $f$ to the difference function $\Delta_{h}f$,
defined by $\displaystyle\Delta_{h}f(x)=f(x+h)-f(x)$. When $h=1$, the
subscript is omitted. This operation is the backbone of the calculus
of finite differences and is fundamental in approximating derivatives,
constructing polynomial interpolants, or
finding closed-form representations of certain sums. In addition, it is an analogue
of the derivative operator of differential calculus.

The factorial powers behave nicely with respect to $\Delta$, as indicated by
the formulas:
\begin{align*}
\Delta x^{\underline{m}}  &  =(x+1)^{\underline{m}}-x^{\underline{m}%
}=(x+1)x^{\underline{m-1}}-x^{\underline{m-1}}(x-m+1)=mx^{\underline{m-1}},\\
\Delta x^{\overline{m}}  &  =(x+1)^{\overline{m}}-x^{\overline{m}%
}=(x+1)^{\overline{m-1}}(x+m)-x(x+1)^{\overline{m-1}}=m(x+1)^{\overline{m-1}}.
\end{align*}
In contrast, the usual powers are not so well suited for finite
differentiation, since%
\[
\Delta x^{m}=(x+1)^{m}-x^{m}=\sum_{k=0}^{m-1}\dbinom{m}{k}x^{k}.
\]
Furthermore, the falling factorial can be \textquotedblleft
normalised\textquotedblright\ in the form $\binom{x}{m}=\frac{x^{\underline{m}%
}}{m!}$, which reproduces the usual binomial coefficients when $x$ is a
positive integer, and satisfies $\Delta\binom{x}{m}=\binom{x}{m-1}$, for
$m\geq1$. This behaviour explains the preference for factorial powers in the
context of successive applications of the difference operator. In this
direction, the following equivalent version of Stirling's formula
(\ref{yanon_eq:02})%
\begin{equation}
x^{m}=\sum_{n=0}^{m}{n}^{\{m\}}\binom{x}{n} \label{yanon_eq:13}%
\end{equation}
is used for expressing polynomials in terms of normalised\ falling factorials,
with subpowers as coefficients. The identity (\ref{yanon_eq:13}) between the two
polynomials of degree $m$ is the consequence of their equality for all
positive integers $x$, as expressed in (\ref{yanon_eq:07}).

The application of $\Delta_{h}$ $n$ times is usually denoted by $\Delta
_{h}^{n}$ (with $\Delta_{h}^{0}$ being the identity operator) and has the
expression%
\begin{equation}
\Delta_{h}^{n}f(x)=\sum_{k=0}^{n}(-1)^{n-k}\binom{n}{k}f(x+kh), \label{yanon_eq:14}%
\end{equation}
which can be proved inductively based on the recurrence relation satisfied by
the binomial coefficients. In light of Proposition \ref{yanon_binom_inv}, $\Delta
_{h}^{n}f(x)$ is the inverse binomial transform of $f(x+nh)$, both understood
as sequences with index $n$, and for fixed values of $x$ and $h$. In
consequence, $f(x+nh)$ is the binomial transform of $\Delta_{h}^{n}f(x)$,
i.e.,%
\begin{equation}
f(x+nh)=\sum_{k=0}^{n}\binom{n}{k}\Delta_{h}^{k}f(x). \label{yanon_eq:15}%
\end{equation}

We recognise in (\ref{yanon_eq:14}) the form of the subpower expression
(\ref{yanon_eq:01}), with $f(x)=x^{m}$, $h=1$ and $x=0$, which means that%
\begin{equation}
n^{\{m\}}=\Delta^{n}x^{m}(0). \label{yanon_eq:16}%
\end{equation}
For this reason, these numbers are referred to in the calculus of finite
differences as the \emph{differences of zero} and are denoted by $\Delta
^{n}0^{m}$.

In \cite{Euler1755}, Euler expressed the formula for the forward difference of higher order
of the monomial $x^{m}$, which, in light of our notation, can be written as%
\begin{equation}
\Delta_{h}^{n}x^{m}=\sum_{k=n}^{m}\dbinom{m}{k}n^{\{k\}}h^{k}x^{m-k}.
\label{yanon_eq:17}%
\end{equation}
This result has a striking resemblance to the binomial expansion of
$(x+nh)^{m}$, which shouldn't be a surprise.
Indeed, based on (\ref{yanon_eq:15}),
$(x+nh)^{m}$ is the binomial transform of $\Delta_{h}^{n}x^{m}$, where the
index is $n$. Writing $(x+nh)^{m}=$ $\sum_{k=0}^{m}\binom{m}{k}n^{k}%
h^{k}x^{m-k}$ and using that $n^{k}$ is the binomial transform of $n^{\{k\}}$,
for every fixed $k$, immediately justifies (\ref{yanon_eq:17}).

\subsection{Some polynomial identities and their ramifications}

The identity (\ref{yanon_eq:13}) could be considered as an alternative definition
of the subpowers. It is also the source of some interesting results.

Look again at the values in Table 1 and perform an \textquotedblleft
alternating right-to-left row sum\textquotedblright. This gives for the first
few lines%
\begin{align*}
2-1  &  =1\\
6-6+1  &  =1\\
24-36+14-1  &  =1\\
120-240+150-30+1  &  =1.
\end{align*}
A clear pattern emerges, which at first looks like a rabbit pulled out of the
hat. The magic trick becomes clear when we replace $x$ by $-x$ in
(\ref{yanon_eq:13}). A simple computation gives $\binom{-x}{n}=(-1)^{n}\binom
{x+n-1}{n}$, which then leads to an equivalent form of (\ref{yanon_eq:13}):
\begin{equation}
x^{m}=\sum_{n=0}^{m}(-1)^{m-n}n^{\{m\}}\binom{x+n-1}{n}. \label{yanon_eq:18}%
\end{equation}
Letting $x=1$ gives the alternate row sum in the array of subpowers%
\[
\sum_{n=0}^{m}(-1)^{m-n}n^{\{m\}}=1,
\]
as predicted.

Now, let us return to (\ref{yanon_eq:18}), written for the exponent $m+1$. Using the
identity $\binom{x+n-1}{n}=\frac{x}{n}\binom{x+n-1}{n-1}$, then shifting the
summation index $n$ by $1$, we obtain
\[
x^{m+1}=\sum_{n=1}^{m+1}(-1)^{m+1-n}n^{\{m+1\}}\frac{x}{n}\binom{x+n-1}{n-1}=x\sum_{n=0}^{m}(-1)^{m-n}\frac{(n+1)^{\{m+1\}}}{n+1}\binom{x+n}{n},
\]
hence
\begin{equation}
x^{m}=\sum_{n=0}^{m}(-1)^{m-n}\dfrac{(n+1)^{\{m+1\}}}{n+1}\binom{x+n}%
{n}\text{.} \label{yanon_eq:19}%
\end{equation}

The signless coefficients $\frac{(n+1)^{\{m+1\}}}{n+1}$ in (\ref{yanon_eq:19}),
expressed in terms of the Stirling set numbers as $n!\,\genfrac{\{}{\}}{0pt}{}{m+1}{n+1}$, were
considered by Worpitzky in \cite{Worpitzky1883} in his study of the Bernoulli
numbers and were referred in \cite{Vandervelde2018} as
\emph{the Worpitzky numbers}, under the notation $W_{m,n}$. Based on
(\ref{yanon_eq:11}) and (\ref{yanon_eq:16}),
\begin{align*}
W_{m,n}  &  =\dfrac{(n+1)^{\{m+1\}}}{n+1}=n^{\{m\}}+(n+1)^{\{m\}}=\Delta^{n}x^{m}(0)+\Delta^{n+1}x^{m}(0)\\
&  =\Delta^{n}\left(  x^{m}+\Delta
x^{m}\right)  (0)=\Delta^{n}\left(  x+1\right)  ^{m}(0)=\Delta^{n}x^{m}(1)=\Delta^{n}1^{m}.
\end{align*}
For this reason, these numbers could also be called \emph{the differences of
one}.

\subsection{Sums of powers and subpowers}

The computation of sums that involve integer powers has been of great interest
to mathematicians since early times. In particular, the sum $S_{m}%
(n)=1^{m}+2^{m}+\ldots+n^{m}$ attracted a lot of attention since antiquity.
Early days formulas\footnote{Without the support of a developed system of
mathematical notations, these \textquotedblleft formulas\textquotedblright%
\ were described using the natural language.} were known for special cases,
when $m$ is $1$, $2$ or $3$, but general results appeared only later, during
the seventeenth and eighteenth centuries. One such general formula expresses $S_{m}(n)$ as
a linear combination of the binomial numbers $\binom{n+1}{1},\binom{n+1}%
{2},\ldots,\binom{n+1}{m+1}$. The idea is to write each power in terms of its
subpowers and then add them up. For $m,n\geq1$, making use of (\ref{yanon_eq:07}) and
Remark \ref{rem:sumation_range} leads to%
\begin{equation}
\sum_{k=1}^{n}k^{m}=\sum_{k=1}^{n}\sum_{p=1}^{m}\binom{k}{p}p^{\{m\}}%
=\sum_{p=1}^{m}\left(  p^{\{m\}}\sum_{k=p}^{n}\binom{k}{p}\right)  =\sum
_{p=1}^{m}p^{\{m\}}\binom{n+1}{p+1} \label{yanon_eq:20}%
\end{equation}
based on the well-known identity $\sum_{k=p}^{n}\binom{k}{p}=\binom{n+1}{p+1}%
$. To exemplify, we can use the values from Table 1 to obtain a formula for
the sum of the fifth powers:%
\begin{equation}
1^{5}+2^{5}+\ldots+n^{5}=\tbinom{n+1}{2}+30\tbinom{n+1}{3}+150\tbinom
{n+1}{4}+240\tbinom{n+1}{5}+120\tbinom{n+1}{6}. \label{yanon_eq:21}%
\end{equation}

Formula (\ref{yanon_eq:20}) is a well-known result and likely the easiest to derive.
However, it is most commonly expressed in terms of the Stirling set numbers
multiplied by factorials. Nevertheless, the form presented here is more appealing
and meaningful, as it connects sums of two different kinds of powers.

A more desirable (and harder to obtain) closed form expression for $S_{m}(n)$
uses powers of $n$, rather than binomial coefficients. For example,
(\ref{yanon_eq:21}) is equivalent to%
\[
1^{5}+2^{5}+\ldots+n^{5}=\frac{1}{6}n^{6}+\frac{1}{2}n^{5}+\frac{5}{12}%
n^{4}-\frac{1}{12}n^{2}\text{.}%
\]
A general formula of this type was obtained independently by the Japanese
mathematician Takakazu Seki~\cite{Seki1712} (1712) and by Jacob Bernoulli~\cite{Bernoulli1713} (1713) (we refer to
\cite{Arakawa2014} for a more detailed account). They noticed the existence of
a sequence of rational numbers $B_{0},B_{1},B_{2},\ldots$ which describe the
coefficients in the expression of $S_{m}(n)$ in a simple uniform manner,
through the formula%
\begin{equation}
S_{m}(n)=\dfrac{1}{m+1}\sum_{k=0}^{m}\binom{m+1}{k}B_{k}n^{m+1-k}.
\label{yanon_eq:22}%
\end{equation}
The coefficients $B_{k}$ in (\ref{yanon_eq:22}) are known today as \emph{the
Bernoulli numbers} (sometimes called \emph{the Seki--Bernoulli numbers}) and
they are defined by the implicit recurrence
\begin{equation}
\sum_{k=0}^{m}\binom{m+1}{k}B_{k}=m+1\qquad(m\geq0), \label{yanon_eq:23}%
\end{equation}
which is just $S_{m}(1)=1$. Solving (\ref{yanon_eq:23}) step by step gives the
values in Table 2. We note the obvious pattern, that
$B_{m}=0$ for odd $m\geq3$, which is a valid general property.

\begin{table}[h]
\caption{Bernoulli numbers.}%
\setlength{\tabcolsep}{7pt}
\renewcommand{\arraystretch}{1.25}
\centering
\par%
\begin{tabular}
[c]{l|ccccccccccccc}%
$m$ & $0$ & $1$ & $2$ & $3$ & $4$ & $5$ & $6$ & $7$ & $8$ & $9$ & $10$ & $11$
& $12$\\\hline
$B_{m}$ & $1$ & $\frac{1}{2}$ & $\frac{1}{6}$ & $0$ & $-\frac{1}{30}$ & $0$ &
$\frac{1}{42}$ & $0$ & $-\frac{1}{30}$ & $0$ & $\frac{5}{66}$ & $0$ &
$-\frac{691}{2730}$%
\end{tabular}
\end{table}

Since their discovery, the Bernoulli numbers have helped to express other
important results, such as the series expansions of some trigonometric and
hyperbolic trigonometric functions, the Euler--Maclaurin summation formula or
the evaluation of the Riemann zeta function, just to name the most common occurrences.
The connection with the Stirling set numbers and, consequently, with subpowers
is established in the following explicit form of the Bernoulli numbers,
written with the help of the subpower notation as%
\begin{equation}
B_{m}=\sum_{n=0}^{m}\frac{(-1)^{m-n}n^{\{m\}}}{n+1} \label{yanon_eq:24}%
\end{equation}
(we refer to \cite[p. 35]{Arakawa2014} for proof and other details). Using the
values from Table 1 we may compute, for example,%
\[
B_{4}=0^{\{4\}}-\dfrac{1^{\{4\}}}{2}+\dfrac{2^{\{4\}}}{3}-\dfrac{3^{\{4\}}}%
{4}+\frac{4^{\{4\}}}{5}=-\frac{1}{2}+\dfrac{14}{3}-\frac{36}{4}+\frac{24}%
{5}=-\dfrac{1}{30}.
\]
At the same time, it is surprising that a very similar sum has a much simpler
answer:%
\begin{equation}
\sum_{n=1}^{m}\dfrac{(-1)^{n+1}n^{\{m\}}}{n}=\delta_{1,m}. \label{yanon_eq:25}%
\end{equation}
This is a direct consequence of (\ref{yanon_eq:19}), written for $x=0$ and with $m$
replaced by $m-1$.

An important variant of the sum $S_{m}(n)$ is the row sum in the array of
subpowers $F(m)=\sum_{n=0}^{m}n^{\{m\}}$. Here we encounter the sequence of
\emph{the Fubini numbers}. The name was proposed by Comtet \cite[p.
238]{Comtet1974} because $F(m)$ represents the number of \textquotedblleft
Fubini iterated formulas\textquotedblright\ for an $m$-dimensional multiple
integral (or sum) that follow from Fubini's theorem. An equivalent formulation
of this combinatorial description is the number of weak orders on $m$ labelled
elements, or, in more familiar terms, the number of ways $m$ competitors can
rank in a competition, where ties are allowed. There is also the analogy with
the row sums in the array of the Stirling set numbers. These sums count the
number of ways to partition a set of a given size and are called \emph{the
Bell numbers}. For this reason, the Fubini numbers are also referred to as
\emph{the ordered Bell numbers}, since they count the number of ways to partition
a set and order the blocks of the partition.

The sequence $\left(  F(m)\right)  _{m\geq0}$ has the entry A000670 in the
OEIS \cite{OEIS} and is governed by the recurrence%
\begin{equation}
F(m+1)=\sum_{k=0}^{m}\binom{m+1}{k}F(k)\qquad(m\geq0),\qquad F(0)=1,
\label{yanon_eq:26}%
\end{equation}
which has a striking resemblance to the corresponding recurrence for the
Bernoulli numbers --- just erase the leftmost $F$ and the following pair of
parentheses and here is (\ref{yanon_eq:23}). To obtain (\ref{yanon_eq:26}), it is enough
to count the number of possible rankings among $m+1$ competitors in a contest
where ties are allowed, by conditioning on the number $k\in\{0,1,\ldots,m\}$
of contestants that \emph{do not} win: they can be chosen in $\binom{m+1}{k}$
ways, and there are $F(k)$ possible rankings among them, while the remaining
$m+1-k$ end up on the first place.

The Fubini numbers arise naturally in combinatorial problems that deal with
set arrangements, orderings, classifications, and hierarchies, but they also
express the results of certain sums or series, like in%
\[
F(m)=\dfrac{1}{2}\sum_{n=0}^{\infty}\dfrac{n^{m}}{2^{n}}\qquad(m\geq0).
\]
This is a particular instance of a more general result%
\begin{equation}
\sum_{n=0}^{\infty}n^{m}x^{n}=\dfrac{1}{1-x}\sum_{k=0}^{m}k^{\{m\}}\left(
\dfrac{x}{1-x}\right)  ^{k},\qquad x\in(0,1),~m\geq0,\label{yanon_eq:27}%
\end{equation}
written for $x=\frac{1}{2}$. To prove (\ref{yanon_eq:27}), it is enough to use
(\ref{yanon_eq:07}) in combination with Remark \ref{rem:sumation_range}, then change
the order of summation (since all terms are nonnegative) to obtain:%
\[
\sum_{n=0}^{\infty}n^{m}x^{n}=\sum_{n=0}^{\infty}\sum_{k=0}^{m}\binom
{n}{k}k^{\{m\}}x^{n}=\sum_{k=0}^{m}k^{\{m\}}\sum_{n=0}^{\infty}\binom
{n}{k}x^{n}=\sum_{k=0}^{m}k^{\{m\}}\dfrac{x^{k}}{(1-x)^{k+1}}.
\]

\section{Down the rabbit hole}

The expression (\ref{yanon_eq:01}) is still legitimate when the exponent takes
positive real values, but to extend this further for negative exponents or complex values,
we have to omit the power of $0$. This adjustment in the
formula only affects the value of $n^{\{0\}}$, due to the missing term
$(-1)^{n}0^{0}=(-1)^{n}$. With this slight change, we can define \emph{the
subpower function} of complex exponent $z$ as
\begin{equation}
n^{\{z\}}=\sum_{k=1}^{n}(-1)^{n-k}\binom{n}{k}k^{z}\qquad(z\in\mathbb{C}),
\label{yanon_eq:28}%
\end{equation}
with $0^{\{z\}}=0$ as an empty sum. This is the same as considering $0^{z}$ to
be $0$ for all $z\in\mathbb{C}$, which modifies the established value of
$n^{\{0\}}$ from $\delta_{0,n}$ to
\[
\delta_{0,n}-(-1)^{n}=\left\{
\begin{tabular}
[c]{rl}%
$0\,$, & $\text{if }n=0$\\
$(-1)^{n+1}\,$, & $\text{if }n\geq1$%
\end{tabular}
\right.  =(-1)^{n+1}(1-\delta_{0,n})
\]
in order to obtain an analytic function over $\mathbb{C}$. When writing
$n^{\{0\}}$, it should become clear from the context, or otherwise be
specified explicitly, which of the two values we mean (see Figure 1).

\begin{figure}[ht]
\caption{The graphs of the subpower functions of real argument $x$, with base
$n$ from $1$ to $5$.}%
\centering
\includegraphics{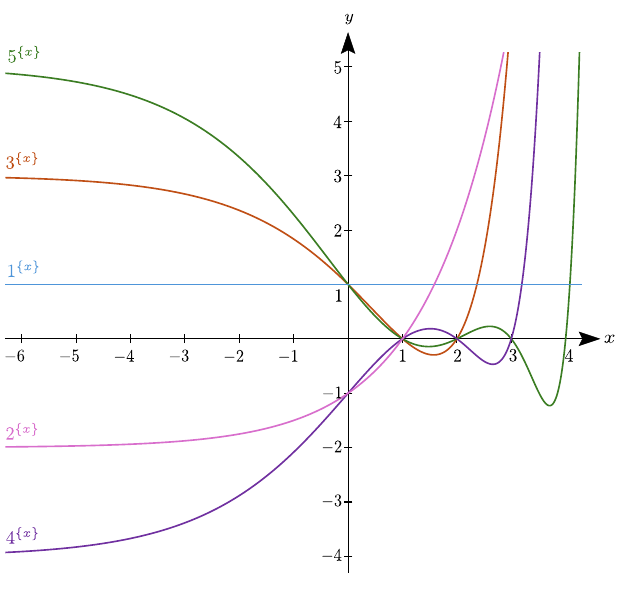}
\end{figure}

Fixing $z$ and agreeing the first term of $\left(
n^{z}\right)  _{n\geq0}$ is $0$, it is still true that its inverse binomial
transform is $\left(  n^{\{z\}}\right)  _{n\geq0}$ with $0^{\{z\}}=0$, which
means that $n^{z}=\sum_{k=1}^{n}\binom{n}{k}k^{\{z\}}$ and $\Delta
^{n}x^{z}(0)=n^{\{z\}}$. Also, the recurrence (\ref{yanon_eq:11}) still
holds for complex exponents:%
\[
n^{\{z\}}=n\left(  n^{\{z-1\}}+(n-1)^{\{z-1\}}\right)  \qquad(n\geq1)\text{,}%
\]
which can be easily checked using the definition (\ref{yanon_eq:28}).
These ideas have already passed the test in the context of the \emph{Stirling
functions}, which are extensions of the Stirling numbers to complex first
argument. These were considered in \cite{Butzer1989,Butzer2003a}, to
which we refer for a detailed account.

\subsection{Subpowers with negative integer exponents}

About three decades ago, Loeb and Rota \cite{Loeb1989} and Roman
\cite{Roman1992,Roman1993} introduced and studied a family of harmonic
numbers to give explicit expressions of the harmonic logarithms. In a
nutshell, when computing the $n$th order antiderivative of $\ln^{m}x$, while
omitting the integration constants, the results can be expressed in a uniform
manner, by%
\begin{equation}
\int\!\!\int\!\!\cdots\!\!\int\ln^{m}x\,(\mathrm{d}x)^{n}=\dfrac{x^{n}}%
{n!}\sum_{k=0}^{m}(-1)^{k}c_{n}^{(k)}m^{\underline{k}}\,\ln^{m-k}x\text{.}
\label{yanon_eq:29}%
\end{equation}
For example,
\[
\iiint\ln^{m}x\,(\mathrm{d}x)^{3}=\frac{x^{3}}{6}\Bigl(\ln^{m}x-\frac{11}{6}m\ln^{m-1}x+\frac{85}{36}m(m-1)\ln^{m-2}x-\frac{575}{216}m(m-1)(m-2)\ln^{m-3}x+\ldots\Bigr).
\]

It turns out that the unspecified coefficients in (\ref{yanon_eq:29}) have the explicit formula%
\[
c_{n}^{(m)}=\sum_{k=1}^{n}\binom{n}{k}\frac{(-1)^{k-1}}{k^{m}}\qquad
(m,n\geq1)
\]
due to Knuth (see \cite[Proposition 3.3.6]{Loeb1989}), hence they can be
expressed as
\begin{equation}
c_{n}^{(m)}=(-1)^{n-1}n^{\{-m\}}.\label{yanon_eq:30}%
\end{equation}
In particular, $c_{n}^{(1)}=1+\frac{1}{2}+\ldots+\frac{1}{n}=H_{n}$ (the $n$th
harmonic number), which means $n^{\{-1\}}=$ $(-1)^{n-1}H_{n}$.

The connections between the harmonic numbers $c_{n}^{(m)}$ and the subpowers
$n^{\{-m\}}$ of negative integer exponent can be exploited by means of the
Stirling set numbers $\genfrac{\{}{\}}{0pt}{}{-m}{n}$ with negative integers in the first
argument, which Branson studied in \cite{Branson1996,Branson2006}. An investigation into
the matter reveals the following equivalent analytic representations:
\begin{align*}
n^{\{-m\}}  &  =\dfrac{(-1)^{m+n-1}\,n!}{m!}\dfrac{\mathrm{d}^{m}}%
{\mathrm{d}x^{m}}\left.  \left(  \frac{1}{(x+1)(x+2)\cdots(x+n)}\right)
\right\vert _{x=0}=\dfrac{(-1)^{n-1}n}{m!}\int_{0}^{1}\!(1-x)^{n-1}\left(  \ln\frac{1}%
{x}\right)  ^{m}\mathrm{d}x\\
&  =\dfrac{(-1)^{n-1}n}{m!}\int_{0}^{1}\!x^{n-1}\!\left(  \ln\frac{1}%
{1-x}\right)  ^{m}\mathrm{d}x=\dfrac{(-1)^{n-1}}{m!}\int_{0}^{1}\!\left(
\ln\frac{1}{1-\sqrt[n]{x}}\right)  ^{m}\mathrm{d}x.
\end{align*}
For proofs and other relevant explanations that go beyond the main purpose of
this paper, we refer to \cite{Loeb1989,Roman1992,Roman1993,Branson1996,Branson2006,Sesma2017}.

\section{Instead of conclusions}

We started small, in search for a notation that best expressed the answer to a classic combinatorial problem.
In the end, we uncovered a web of connections with various areas of mathematics.
So, what better closing argument than a final connection to a famous theorem? 
The reader is thus invited to reflect on the following version of \emph{Fermat's Last Theorem}, stated for subpowers.

\begin{theorem}
[\footnote{There exists a truly marvellous proof of this, which this footnote is too short to contain.}]
For any positive integer $m\notin\{2,5,7\}$, the equation%
\[
x^{\{m\}}+y^{\{m\}}=z^{\{m\}}%
\]
has no solutions\footnote{$(1,1,2)$ is the solution when $m=2$. If $m=5$, then $(2,5,3)$ and $(5,5,4)$ are solutions. For $m=7$, the solution is $(4,4,5)$.} $x,y,z\in\{1,2,\ldots,m\}$.
\end{theorem}

{\small
\bibliographystyle{habbrv}
\bibliography{yanon.bib}}

\end{document}